\documentclass[12pt,oneside,english]{amsart}
\usepackage[T1]{fontenc}
\usepackage[utf8]{inputenc}
\usepackage[a4paper]{geometry}
\usepackage{verbatim}
\usepackage{mathtools}
\usepackage{amstext}
\usepackage{amsthm}
\usepackage{amssymb}
\usepackage{stmaryrd}

\makeatletter
\numberwithin{equation}{section}
\numberwithin{figure}{section}
\theoremstyle{plain}
\newtheorem{thm}{\protect\theoremname}
  \theoremstyle{plain}
  \newtheorem{prop}[thm]{\protect\propositionname}
  \theoremstyle{remark}
  \newtheorem{rem}[thm]{\protect\remarkname}
  \theoremstyle{plain}
  \newtheorem{cor}[thm]{\protect\corollaryname}
  \theoremstyle{plain}
  \newtheorem{lem}[thm]{\protect\lemmaname}
  \theoremstyle{definition}
  \newtheorem{defn}[thm]{\protect\definitionname}

\usepackage{amstext}

\numberwithin{equation}{section}
\numberwithin{figure}{section}

\usepackage{amstext}

\numberwithin{equation}{section}
\numberwithin{figure}{section}
\theoremstyle{plain}

\usepackage{amstext}

\numberwithin{equation}{section}
\numberwithin{figure}{section}

\usepackage{amstext}

\usepackage{aecompl}

\@namedef{subjclassname@2020}{\textup{2020} Mathematics Subject Classification}

\usepackage{amstext}
\usepackage[all]{xy}

\def\quot{/\!\!/}

\def\sym{\mathsf{Sym}}
\def\hilb{\mathsf{Hilb}}

\def\hom{\mathsf{Hom}}
\def\codim{\mathsf{codim}}

\def\Pexp{\mathsf{PExp}}
\def\Plog{\mathsf{PLog}}

\def\GL{\mathrm{GL}_n}

\def\HE{\mathsf{HE}}

\usepackage{babel}

\title[Moduli spaces of Higgs bundles over abelian varieties]
{Topology of the moduli spaces of \\ Higgs bundles over abelian varieties}

\author[I. Biswas]{Indranil Biswas}

\address{Mathematics Department, Shiv Nadar University, NH91, Tehsil
Dadri, Greater Noida, Uttar Pradesh 201314, India}

\email{indranil29@gmail.com, indranil@math.tifr.res.in}

\author[C. Florentino]{Carlos Florentino}

\address{Departamento de Matem\'{a}tica, Faculdade de Ci\^{e}ncias, Univ. 
de Lisboa, Edf. C6, Campo Grande 1749-016 Lisboa, Portugal}

\email{caflorentino@fc.ul.pt}

\author[A. Nozad]{Azizeh Nozad}

\address{School of Mathematics, Institute for Research in Fundamental Sciences (IPM), P.O. Box 19395-5746, Tehran, Iran}

\email{anozad@ipm.ir}



\newcommand{\xdownarrow}[1]{%
  {\left\downarrow\vbox to #1{}\right.\kern-\nulldelimiterspace}
}

\usepackage{babel}

\makeatother

\usepackage{babel}
  \providecommand{\corollaryname}{Corollary}
  \providecommand{\definitionname}{Definition}
  \providecommand{\lemmaname}{Lemma}
  \providecommand{\propositionname}{Proposition}
  \providecommand{\remarkname}{Remark}
\providecommand{\theoremname}{Theorem}

\begin{document}
\begin{abstract}
Let $G$ be a complex reductive group and $A$ be an Abelian variety
of dimension $d$ over $\mathbb{C}$. We determine the Poincaré polynomial
and also the mixed Hodge polynomial of the moduli space $\mathcal{M}_{A}^{H}(G)$
of $G$-Higgs bundles over $A$. We show that these are normal varieties
with symplectic singularities, when $G$ is a classical semisimple
group. For $G\,=\,\text{GL}_{n}(\mathbb{C})$, we also compute Poincaré
polynomials of natural desingularizations of $\mathcal{M}_{A}^{H}(G)$
and of $G$-character varieties of free abelian groups, in some cases.
In particular, explicit formulas are obtained when $\dim A\,=\,d\,=\,1$,
and also for rank 2 and 3 Higgs bundles, for arbitrary $d\,>\,1$. 
\end{abstract}

\dedicatory{Dedicated to Peter E. Newstead}
\maketitle

\section{Introduction}

The study of the topology of moduli spaces of Higgs bundles over algebraic
curves was initiated in the very first article on Higgs bundles, \cite{Hi},
where these geometric objects were defined by Hitchin. In particular,
in that work and in the work of Gothen, \cite{Got}, we find some
of the first computations of their Poincaré polynomials. Most of the
stronger results deal with the cases where these moduli spaces are
smooth, and the curve has genus $>1$.

More recently, the study of moduli spaces of $G$-Higgs bundles over
elliptic curves, for a complex reductive group $G$, was started by
Franco, Garcia-Prada and Newstead in \cite{FGPN}, and has been generalized
to $G$-Higgs bundles over abelian varieties by Biswas \cite{BF,Bi}
and Franco-Tortella \cite{FT}.

In this context, the non-abelian Hodge correspondence — between the
moduli spaces of $G$-Higgs bundles and the $G$-character varieties
of surface groups — became a very important tool in understanding
the geometry of both these classes of spaces. The study of character
varieties of \emph{free abelian} groups was initiated by Thaddeus
\cite{Th} and Sikora \cite{Sk}, from the point of view of the algebraic
geometry and relation to mirror symmetry, and by Florentino, Lawton
and Silva in \cite{FL,FS} from the point of view of the underlying
topology and mixed Hodge structures. These character varieties are
generally singular but of orbifold type.

Take a complex abelian variety $A$ of dimension $d$ and a reductive
complex affine algebraic group $G$. For a semistable $G$-Higgs bundles
over $A$, the underlying principal $G$-bundle semistable. Let $\mathcal{M}_{A}^{H}(G)$
denote the moduli space of semistable $G$-Higgs bundles over $A$
whose underlying principal $G$-bundle is in the connected component,
containing the trivial principal $G$-bundle, of the moduli space
of principal $G$-bundles. In this article we prove that, when $G$
is $\GL(\mathbb{C})$ or a classical semisimple group, the moduli
spaces $\mathcal{M}_{A}^{H}(G)$ have symplectic singularities as
defined by Beauville in \cite{Be2}. We also determine the Poincaré
polynomials, as well as mixed Hodge polynomials, of $\mathcal{M}_{A}^{H}(G)$
and of the corresponding $G$-character varieties, namely the moduli
spaces of equivalence classes of $G$-representations of $\pi_{1}(A)\,\cong\,\mathbb{Z}^{2d}$.

In the last section, we consider $G=\GL(\mathbb{C})$ and prove that
$\mathcal{M}_{A}^{H}(G)$ (as well as the $\GL(\mathbb{C})$-character
varieties of the free abelian groups) admit a symplectic (also called
crepant) resolution only when $d=1$, and that in this case the resolution
is given by the Hilbert-Chow morphism.

We finish with an explicit calculation of the Poincaré polynomials
of the resolutions of $\mathcal{M}_{A}^{H}(G)$ and of the $G$-character
varieties of $\mathbb{Z}^{2d}$, in all cases when their resolutions
are Hilbert schemes. These are precisely the crepant cases: $d=1$
and arbitrary rank $n\in\mathbb{N}$; and the non-crepant cases: $n=2$
or $n=3$ and arbitrary $d>1$.

\subsection*{Acknowledgements}

We are grateful to the anonymous referee for detailed comments and
for suggesting important corrections on a previous version. We thank
Giosuè Muratore, Jaime Silva, Alfonso Zamora for useful discussion
on these topics. IB is partially supported by a J. C. Bose Fellowship.
CF is partially supported by FCT project UIDB/04561/2020. AN is supported
by a grant from IPM, Iran.

\section{Preliminaries on Character Varieties and Higgs Bundles}

As before, $G$ is a complex reductive affine algebraic group. Its
Lie algebra will be denoted by $\mathfrak{g}$. Let $X$ be an irreducible
smooth projective variety over $\mathbb{C}$, with fundamental group
$\pi_{1}(X)$; we will omit reference to the base point.

\subsection{Character varieties}

Since $\pi_{1}(X)$ is a finitely presented group, the space $\hom(\pi_{1}(X),\,G)$
of group homomorphisms $\pi_{1}(X)\,\longrightarrow\,G$ has a natural
structure of an affine complex algebraic variety. The adjoint action
of $G$ on itself produces an action of $G$ on $\hom(\pi_{1}(X),\,G)$.
The (affine) geometric invariant theoretic quotient 
\[
\hom(\pi_{1}(X),\,G)\quot G,
\]
which is called the $G$-character variety of $\pi_{1}(X)$, is an
affine variety whose algebra of regular functions consists of the
$G$-invariant regular functions on $\hom(\pi_{1}(X),G)$.

Let $A$ be a complex abelian variety of dimension $d$. Denote by
\begin{equation}
\mathcal{M}_{A}^{R}(G)\,:=\,\hom(\pi_{1}(A),G)^{0}\quot G\label{eq:char-var}
\end{equation}
the connected component, containing the trivial homomorphism, of the
$G$-character variety of $\pi_{1}(A)\cong\mathbb{Z}^{2d}$. It is
an irreducible affine variety of dimension $2dr$, where $r$ is the
rank of $G$, as shown in \cite[Thm. 2.1(1)]{Sk} and \cite[Cor. 5.14]{FL}.

The celebrated non-abelian Hodge correspondence produces a $C^{\infty}$
identification between the $G$-character variety $\hom(\pi_{1}(X),G)\quot G$
and a moduli space of $G$-Higgs bundles over $X$, whose definition
we now recall.

\subsection{\label{subsec:PHB}Principal Higgs bundles over abelian varieties}

For any principal $G$-bundle $E$ on a smooth complex projective
variety $X$, the adjoint bundle $E(\mathfrak{g})$ is the Lie algebra
bundle on $X$ associated to $E$ for the $G$-module $\mathfrak{g}\,=\,\text{Lie}(G)$
defined by the adjoint action. The vector bundle over $X$ whose sheaf
of sections is given by the holomorphic $i$-forms will be denoted
by $\Omega_{X}^{i}$.

Recall that a $G$-Higgs bundle over $X$ is a pair $(E,\,\theta)$
where $E$ is a holomorphic principal $G$-bundle over $X$ and $\theta$
is a holomorphic section of $E(\mathfrak{g})\otimes\Omega_{X}^{1}$
satisfying the equation $\eta(\theta\otimes\theta)\,=\,0$, where
\[
\eta:(E(\mathfrak{g})\otimes\Omega_{X}^{1})^{\otimes2}\,\longrightarrow\,E(\mathfrak{g})\otimes\Omega_{X}^{2}
\]
is the homomorphism obtained by combining the Lie bracket operation
on $E(\mathfrak{g})$ with the exterior product of 1-forms.

Fix an ample line bundle on $X$ to define the degree of coherent
sheaves on $X$. See \cite{BG}, \cite{Si0}, \cite{Sc} for the definition
of semistable $G$-Higgs bundles. Let $\mathcal{\widetilde{M}}_{X}^{H}(G)$
denote the moduli space of semistable $G$-Higgs bundles with vanishing
rational characteristic classes \cite{Sc}, \cite{Si}.

There are homeomorphisms 
\begin{equation}
\hom(\pi_{1}(X),G)\quot G\,\,\cong\,\,\,\mathcal{\widetilde{M}}_{X}^{H}(G),\label{eq:NAHC}
\end{equation}
given by the non-abelian Hodge theory (see \cite{Si0,Hi,BG}). In
particular, these two spaces have the same homotopy type as well as
the rational cohomologies.

We will now consider $G$-Higgs bundles over abelian varieties.

As before, let $A$ be a complex abelian variety of dimension $d$.
For any semistable $G$-Higgs bundle $(E,\,\theta)$ on $A$, the
underlying principal $G$-bundle $E$ is semistable \cite{BG,BF}.
Denote by $\mathcal{M}_{A}(G)$ the connected component, containing
the trivial principal $G$-bundle, of the moduli space of semistable
topologically trivial principal $G$-bundles on $A$. Also, denote
by $\mathcal{M}_{A}^{H}(G)$ the moduli space of semistable $G$-Higgs
bundles over $A$ such that the underlying principal $G$-bundle is
in $\mathcal{M}_{A}(G)$.

\section{Normality of the Moduli Spaces}

Let $G$ be a connected reductive group over $\mathbb{C}$. We start
by recalling the descriptions of $\mathcal{M}_{A}^{H}(G)$ and of
$\mathcal{M}_{A}^{R}(G)$ in \cite{Bi}.

Let $A^{\vee}\,:=\,{\rm Pic}^{0}(A)$ be the dual abelian variety,
parametrizing the isomorphism classes of topologically trivial holomorphic
line bundles on $A$. Fix a polarization $\xi\,\in\,H^{2}(A;\,\mathbb{Q})\bigcap H^{1,1}(A)$.
So we have a nondegenerate pairing 
\[
H^{0}(A,\,\Omega_{A}^{1})\,\otimes\,H^{1}(A,\,\mathcal{O}_{A})\,\longrightarrow\,\mathbb{C},\,\ \ \alpha\otimes\beta\,\longmapsto\int_{A}\alpha\wedge\beta\wedge\xi^{d-1}.
\]
This provides identifications 
\begin{equation}
H^{0}(A,\,\Omega_{A}^{1})\,=\,H^{1}(A,\,\mathcal{O}_{A})^{*}\,=\,T_{z}^{*}A^{\vee}\label{td}
\end{equation}
for all $z\,\in\,A^{\vee}$.

Now, fix a maximal torus $T\,\subset\,G$ with the corresponding Weyl
group $W$, and denote the multiplicative group $\mathbb{C}^{*}$
by $\mathbb{G}_{m}$. Then, the \emph{cocharacter lattice} of $T$
\[
\Lambda_{T}:=\hom(\mathbb{G}_{m},\,T)
\]
is naturally a $W$-module. Abbreviate by $\mathcal{M}^{H}:=\mathcal{M}_{A}^{H}(\mathbb{G}_{m})$
the moduli space of topologically trivial $\mathbb{G}_{m}$-Higgs
bundles over $A$. From \eqref{td} we have $\mathcal{M}^{H}\,\cong\,T^{*}A^{\vee}$,
the cotangent bundle of $A^{\vee}$. Recall that $\mathcal{M}_{A}^{R}(G)$
is the character variety $\hom(\mathbb{Z}^{2d},G)^{0}\quot G$ (as
in \eqref{eq:char-var}), and let us denote $\mathcal{M}_{A}^{R}(\mathbb{G}_{m})\cong\mathbb{G}_{m}^{2d}$
by $\mathcal{M}^{R}$. 
\begin{thm}[{\cite{Bi}}]
\label{thm:M^H} There are natural algebraic bijections 
\[
\varphi\,:\,(\mathcal{M}^{H}\otimes_{\mathbb{Z}}\Lambda_{T})/W\,\longrightarrow\,\mathcal{M}_{A}^{H}(G),
\]
and 
\[
\psi\,:\,(\mathcal{M}^{R}\otimes_{\mathbb{Z}}\Lambda_{T})/W\,\longrightarrow\,\mathcal{M}_{A}^{R}(G),
\]
where the Weyl group $W$ acts naturally on $\Lambda_{T}$ and trivially
on $\mathcal{M}^{H}$ and $\mathcal{M}^{R}$. 
\end{thm}

In the case $G={\rm GL}_{n}(\mathbb{C})$, the moduli spaces in Theorem
\ref{thm:M^H} can be identified with symmetric products of smooth
varieties. Given a quasi-projective complex variety $X$, we denote
by $\sym^{n}(X)$ its $n$th symmetric product, that is the quotient
$X^{n}/S_{n}$ where the symmetric group $S_{n}$ acts on $X^{n}$
by permuting the entries. 
\begin{prop}
\label{prop:GLn}Let $G\,=\,{\rm GL}_{n}(\mathbb{C})$ and $A$ an
abelian variety of complex dimension $d$. Then there are natural
isomorphisms: 
\[
\mathcal{M}_{A}^{H}(G)\cong\sym^{n}(\mathcal{M}^{H}),
\]
and 
\[
\mathcal{M}_{A}^{R}(G)\cong\sym^{n}(\mathbb{G}_{m}^{2d}).
\]
\end{prop}

\begin{proof}
This follows from \cite[Thm. 1.1]{FT}; we provide a sketch of the
arguments for completeness. Let us start with the second isomorphism.
In this case, the character variety is irreducible (see \cite[Thm. 2.1(1)]{Sk}),
and we have 
\begin{equation}
\mathcal{M}_{A}^{R}(G)\,=\,\hom(\mathbb{Z}^{2d},\,G)/G\,=\,T^{2d}/W,\label{eq:1-conn-G}
\end{equation}
where $T\,=\,\mathbb{G}_{m}^{n}$ is a maximal torus of $G\,=\,{\rm GL}_{n}(\mathbb{C})$
and $W$ its Weyl group (see \cite{FS}). Since $W\,=\,S_{n}$, the
symmetric group on $n$ elements, acts diagonally on $T^{2d}$ using
the conjugation action of $W$ on $T$, we get: 
\begin{equation}
T^{2d}/W\,=\,(\mathbb{G}_{m}^{n})^{2d}/S_{n}\,\cong\,\sym{}^{n}(\mathbb{G}_{m}^{2d}).\label{eq:Sym^n}
\end{equation}
A similar construction works also for the moduli space of topologically
trivial ${\rm GL}_{n}(\mathbb{C})$-Higgs bundles over $A$, which
are rank $n$ Higgs vector bundles. Indeed, a semistable rank $n$
Higgs bundle over $A$, with topologically trivial underlying vector
bundle, is $S$-equivalent to a direct sum of $\mathbb{C}^{*}$-Higgs
bundles. The ordering in the direct sum is irrelevant, so we get:
\[
\mathcal{M}_{A}^{H}(G)\,\cong\,(\mathcal{M}^{H})^{n}/S_{n}\,=\,\sym^{n}(\mathcal{M}^{H})
\]
as in the statement of the proposition. 
\end{proof}
\begin{rem}
For $G={\rm GL}_{n}(\mathbb{C})$, the descriptions in Theorem \ref{thm:M^H}
and in Proposition \ref{prop:GLn} are indeed compatible, since the
Weyl group $W\,=\,S_{n}$ acts on the lattice $\Lambda_{T}\,=\,\hom(\mathbb{G}_{m},\,\mathbb{G}_{m}^{n})$,
which is isomorphic to $\mathbb{Z}^{n}$, by permuting the standard
generators. Consequently, $(\mathcal{M}^{H}\otimes_{\mathbb{Z}}\Lambda_{T})/S_{n}\,\cong\,\sym^{n}(\mathcal{M}^{H})$,
and similarly for $\mathcal{M}_{A}^{R}(G)$. 
\end{rem}

Let us say that $G$ is a classical semisimple group if it is one
of ${\rm SL}_{n}(\mathbb{C})$, ${\rm Sp}_{n}(\mathbb{C})$ or ${\rm SO}_{n}(\mathbb{C})$,
for some $n\in\mathbb{N}$. Combining the isosingularity theorem of
Simpson \cite[Thm. 10.6]{Si}, with results of Sikora \cite{Sk},
we can show that the morphisms in Theorem \ref{thm:M^H} are actually
isomorphisms in the classical semisimple case. 
\begin{thm}
\label{thm:M^H2} Let $G$ be a classical semisimple group. Then the
maps $\varphi$ and $\psi$ in Theorem \ref{thm:M^H} are algebraic
isomorphisms. 
\end{thm}

\begin{proof}
We start with the case of $\mathcal{\mathcal{M}}_{A}^{R}(G)$. For
a maximal torus $T\subset G$ with Weyl group $W$, we consider the
composition of maps 
\[
T^{2d}\,=\,\hom(\mathbb{Z}^{2d},\,T)\,\subset\,\hom(\mathbb{Z}^{2d},\,G)\,\twoheadrightarrow\,\mathcal{\mathcal{M}}_{A}^{R}(G)\,=\,\hom(\mathbb{Z}^{2d},\,G)\quot G.
\]
This factors through a morphism $\chi\,:\,T^{2d}/W\,\longrightarrow\,\mathcal{\mathcal{M}}_{A}^{R}(G)$,
where $W$ acts diagonally on $T^{2d}$ using the natural action of
$W$ on $T$. Moreover, $\chi$ is an isomorphism when $G$ is a classical
semisimple group by \cite[Thm 2.1(3)]{Sk}. Since $T^{2d}/W$ is normal,
so is $\mathcal{\mathcal{M}}_{A}^{R}(G)$. Being a bijective morphism
to a normal target, by Zariski's main theorem we obtain that $\psi$
is an isomorphism.

Now, the non-abelian Hodge correspondence between $\mathcal{M}^{H}(G)$
and $\mathcal{M}^{R}(G)$, and the same correspondence for $\mathbb{C}^{*}$,
allows us to build the following diagram 
\begin{eqnarray*}
(\mathcal{M}^{H}\otimes_{\mathbb{Z}}\Lambda_{T})/W & \longrightarrow & (\mathcal{M}^{R}\otimes_{\mathbb{Z}}\Lambda_{T})/W\\
\varphi\big\downarrow\quad &  & \cong\big\downarrow\psi\\
\mathcal{M}_{A}^{H}(G) & \longrightarrow & \mathcal{M}_{A}^{R}(G),
\end{eqnarray*}
whose commutativity can be easily checked. By Simpson's isosingularity
theorem \cite[Thm. 10.6]{Si}, both horizontal arrows are local isomorphisms
for the étale topology (see also \cite[Rem. 3.3]{Ti}). Since a bijective
algebraic morphism which is étale is an isomorphism, we conclude that
$\varphi$ is an algebraic isomorphism. 
\end{proof}
The following is immediate from Proposition \ref{prop:GLn} and Theorem
\ref{thm:M^H2}.
\begin{cor}
\label{cor-nor} Assume that $G\,=\,\GL(\mathbb{C})$ or it is a classical
semisimple group. Then $\mathcal{M}_{A}^{H}(G)$ and $\mathcal{M}_{A}^{R}(G)$
are irreducible and normal varieties. 
\end{cor}

\begin{rem}
This corollary extends \cite[Thm 3.7]{FGPN} to the context of general
abelian varieties. The question of normality of $\mathcal{M}_{A}^{H}(G)$
and $\mathcal{M}_{A}^{R}(G)$ for general reductive $G$ remains open,
as far as we know. On the other hand, \cite[Thm 4.16]{La} showed
that the moduli spaces of principal $G$-bundles $\mathcal{M}_{A}(G)$
are irreducible normal varieties, for general $G$, when $\dim A=1$. 
\end{rem}

\section{Mixed Hodge Polynomials}

For background on mixed Hodge structures (abbreviated as MHS), we
refer to \cite{PS} and \cite{D}. Given a quasi-projective algebraic
variety $X$ over $\mathbb{C}$, its mixed Hodge polynomial (MHP)
$\mu_{X}$ is the polynomial in three variables whose coefficients
are the Deligne-Hodge numbers 
\[
h^{k,p,q}(X)=\dim H^{k,p,q}(X),
\]
where $H^{k,p,q}(X)$ is the $(p,q)$ graded piece of $H^{k}(X;\,\mathbb{Q})$,
the $k$th cohomology of $X$, with rational coefficients: 
\begin{equation}
\mu_{X}(t,u,v)=\sum_{k,p,q\geq0}h^{k,p,q}(X)\,t^{k}\,u^{p}\,v^{q}.\label{eq:mu}
\end{equation}
From $\mu_{X}$ we can obtain the Poincaré polynomial $P_{X}$ and
the Euler characteristic $\chi_{X}$ as specializations: $P_{X}(t)\,:=\,\sum_{k\geq0}b^{k}(X)\,t^{k}=\mu_{X}(t,1,1)$,
where $b^{k}(X)\,=\,\dim H^{k}(X;\,\mathbb{Q})$ is the $k$-th Betti
number. Note that $\chi_{X}\,=\,\mu_{X}(-1,1,1)$.

Let $T$ be a maximal torus of $G$. Its Lie algebra will be denoted
by $\mathfrak{t}$. The Weyl group $W$ for $T$ acts on both $T$
and on $\mathfrak{t}$. The natural morphism from $T^{2d}/W$ to $\mathcal{M}_{A}^{R}(G)$
\cite[Thm. 2.1]{Sk} induces an isomorphism of mixed Hodge structures;
this yields the following result. 
\begin{thm}
$($\cite[Thm. 5.2]{FS}, \cite[Thm. 5.8]{FLS}$)$ Let $G$ be a
connected reductive group. The mixed Hodge polynomial of $\mathcal{M}_{A}^{R}(G)$
is given by: 
\[
\mu_{\mathcal{M}_{A}^{R}(G)}(t,\,u,\,v)\,=\,\frac{1}{|W|}\sum_{w\in W}(\det(I+t\,u\,v\,w))^{2d},
\]
where $w\,\in\,W$ acts on $\mathfrak{t}$ and $I$ is the identity
automorphism of $\mathfrak{t}$. 
\end{thm}

In particular, the Poincaré polynomial of $\mathcal{M}_{A}^{R}(G)$
is (see also \cite{FLS}): 
\begin{equation}
P_{\mathcal{M}_{A}^{R}(G)}(t)=\frac{1}{|W|}\sum_{w\in W}(\det(I+t\,w))^{2d}.\label{eq:PP-R}
\end{equation}

Now, we can state the analogous result for the mixed Hodge polynomials
of $\mathcal{M}_{A}^{H}(G)$ and of $\mathcal{M}_{A}(G)$, the connected
component (of the trivial $G$-bundle) of the moduli space of semistable
topologically trivial $G$-bundles on $A$. 
\begin{thm}
The mixed Hodge structures on $\mathcal{M}_{A}(G)$ and on $\mathcal{M}_{A}^{H}(G)$
coincide, and their mixed Hodge polynomial is given by: 
\[
\mu_{\mathcal{M}_{A}^{H}(G)}(t,\,u,\,v)\,=\,\frac{1}{|W|}\sum_{w\in W}(\det(I+t\,u\,w)\det(I+t\,v\,w))^{d},
\]
where $I$ is the identity automorphism on $\mathfrak{t}$. 
\end{thm}

\begin{proof}
By \cite[Thm. 2.3.5]{D}, any morphism of MHSs is strict with respect
to both the weight and the Hodge filtration. This means that any algebraic
map that induces an isomorphism in usual cohomology, induces also
an isomorphism of MHSs. Now, the algebraic map 
\begin{equation}
\mathcal{M}_{A}(G)\hookrightarrow\mathcal{M}_{A}^{H}(G)\label{eq:inclusion}
\end{equation}
given by the natural inclusion $[E]\,\longmapsto\,[(E,\,0)]$ for
a $G$-bundle $E$, is a strong deformation retraction (see \cite{BF}),
so that these moduli spaces have isomorphic MHSs on their cohomologies.

In the same way, the algebraic bijection \cite[Thm. 2.2]{Bi} 
\[
\delta:(A^{\vee}\otimes_{\mathbb{Z}}\Lambda_{T})/W\,\longrightarrow\,\mathcal{M}_{A}(G)
\]
which induces isomorphisms in cohomology (see also \cite{BFM}) gives
also isomorphisms of MHSs: 
\[
H^{*}(\mathcal{M}_{A}(G))\cong H^{*}\left((A^{\vee}\otimes_{\mathbb{Z}}\Lambda_{T})/W\right)\cong H^{*}(A^{\vee}\otimes_{\mathbb{Z}}\Lambda_{T})^{W},
\]
where the last space denotes the $W$-invariant subspaces of $H^{*}(A^{\vee}\otimes_{\mathbb{Z}}\Lambda_{T}),$
from which $\mu_{\mathcal{M}_{A}(G)}$ can now be obtained.

Observe now that $B:=A^{\vee}\otimes_{\mathbb{Z}}\Lambda_{T}$ is
also an abelian variety of dimension $n=dr$, where $r=\dim T$ is
the rank of $G$. Moreover, it is well known that the mixed Hodge
polynomial of any abelian variety $B$ of dimension $n$ is: 
\[
\mu_{B}(t,\,u,\,v)\,=\,\left((1+tu)(1+tv)\right)^{n},
\]
since $B$ has a pure Hodge structure and it is generated in degree
1; moreover, in the decomposition 
\[
H^{1}(B)\,=\,H^{1}(A^{\vee}\otimes_{\mathbb{Z}}\Lambda_{T})\,=\,H^{1,0}(B)\oplus H^{0,1}(B)
\]
$H^{1,0}(B)$ can be naturally identified with the direct sum of $d\,=\,\dim A$
copies of $\mathfrak{t}$, the Lie algebra of $T$. Now, recall that
if $V=\oplus_{k\geq0}V^{k}$ is a (finite dimensional) graded $\mathbb{C}$-vector
space and $g\,:\,V\,\longrightarrow\,V$ is a linear map preserving
the grading, the graded-character of $g$, defined as the series in
the formal variable $t$: 
\begin{equation}
\chi_{g}(V):=\sum_{k\geq0}\text{tr}(g|_{V^{k}})\,t^{k}\in\mathbb{C}[t]\label{eq:char}
\end{equation}
is additive and multiplicative, under direct sums and tensor products,
respectively, and verifies: 
\[
\chi_{g}(\wedge^{*}V)=\det(I+t^{\delta}g),
\]
whenever all elements of $V$ have degree $\delta\in\mathbb{Z}$ (see,
for example \cite{Se}). In our case, we can write: 
\[
H^{*}(B)\cong\bigwedge^{*}H^{1}(B)\cong\bigwedge^{*}H^{1,0}(B)\otimes\bigwedge^{*}H^{0,1}(B)=\bigwedge^{*}\mathfrak{t}^{d}\otimes\bigwedge^{*}\overline{\mathfrak{t}}^{d}\cong(\bigwedge^{*}\mathfrak{t}\otimes\bigwedge^{*}\overline{\mathfrak{t}})^{\otimes d},
\]
where, in terms of the triple grading, $\mathfrak{t}$ has degree
$(1,1,0)$ and $\overline{\mathfrak{t}}$ has degree $(1,0,1)$ (corresponding,
respectively, to the monomials $tu$ and $tv$ in \eqref{eq:mu}).
Hence, applying \eqref{eq:char} to the computation of the dimension
of the $W$ representation $H^{*}(B)$ (see, for example \cite[Cor. 4.4]{FLS}),
we conclude that 
\[
\mu_{B/W}(t,u,v)=\frac{1}{|W|}\sum_{w\in W}\chi_{w}(H^{*}(B))=\frac{1}{|W|}\sum_{w\in W}\det(I+tu\,w)^{d}\det(I+tv\,w)^{d}
\]
as wanted. 
\end{proof}
\begin{rem}
Note that this formula specializes, for $u=v=1$, to the same Poincaré
polynomial as in Equation \eqref{eq:PP-R}, as it should, given the
homeomorphism $\mathcal{M}_{A}^{H}(G)\simeq\mathcal{M}_{A}^{R}(G)$.

On the other hand, the mixed Hodge structures on $\mathcal{M}_{A}^{H}(G)$
and on $\mathcal{M}_{A}^{R}(G)$ are very different: $\mathcal{M}_{A}^{R}(G)$
is a Stein analytic space (which implies that $\mathcal{M}_{A}^{R}(G)$
has no positive dimensional projective subvarieties), whereas $\mathcal{M}_{A}^{H}(G)$
has projective subvarieties, such as $\mathcal{M}_{A}(G)$. In particular,
the MHS on $\mathcal{M}_{A}^{H}(G)$ is actually pure and, by contrast,
$\mathcal{M}_{A}^{R}(G)$ has what is called a \emph{balanced} or
\emph{Hodge-Tate} MHS. 
\end{rem}

In general, even though $\mathcal{M}^{H}$ (see Theorem \ref{thm:M^H})
is a smooth quasi-projective variety, both $\mathcal{M}_{A}^{R}(G)$
and $\mathcal{M}_{A}^{H}(G)$ are generally singular algebraic varieties
of orbifold type. In some cases, one can find natural desingularizations,
and compute their MHP / Poincaré polynomials. We first consider their
singularity types.

\section{Symplectic singularities and resolutions}

Let $X$ be a normal and irreducible complex algebraic variety. We
say that $X$ is a \emph{symplectic variety} if its smooth part $X^{sm}\subset X$
admits a holomorphic symplectic form $\omega$ (closed and non-degenerate)
whose pull-back to any resolution $\pi\,:\,Y\,\longrightarrow\,X$
extends to a holomorphic 2-form $\widetilde{\omega}$ on $Y$ (not
necessarily non-degenerate). Equivalently, by Namikawa's theorem (\cite[Thm. 6]{Na1}),
$X$ is symplectic if and only if it has only rational Gorenstein
singularities and its smooth part carries a holomorphic symplectic
form.

If $\pi\,:\,Y\,\longrightarrow\,X$ is a resolution of singularities
of a symplectic variety $X$, then $\pi$ is called a \emph{symplectic
resolution} if the above extended form $\widetilde{\omega}$ is non-degenerate
(so that it is a \emph{de facto} symplectic form on $Y$). This notion
is actually independent of the starting symplectic form $\omega$
on $X^{sm}$. Therefore, if $\pi\,:\,Y\,\longrightarrow\,X$ is a
symplectic resolution, then the canonical divisor $K_{Y}$ is trivial.

\subsection{Symplectic singularities}
\begin{thm}
Assume that $\mathcal{M}_{A}^{H}(G)$ is normal. Then, $\mathcal{M}_{A}^{H}(G)$
and $\mathcal{M}_{A}^{R}(G)$ are symplectic varieties, and have only
rational Gorenstein singularities. 
\end{thm}

\begin{proof}
Let us first consider the $G$-Higgs bundle case. When $G={\rm GL}_{n}(\mathbb{C})$
the theorem follows from the fact that symmetric products of symplectic
varieties are symplectic (see \cite[Prop. 2.4]{Be2}). Since $\mathcal{M}_{A}^{H}(G)$
is normal, we have the following isomorphism 
\[
\mathcal{M}_{A}^{H}(G)\,\cong\,(\mathcal{M}^{H}\otimes_{\mathbb{Z}}\Lambda_{T})/W,
\]
from Theorem \ref{thm:M^H2}. It is well known that $\mathcal{M}^{H}$,
the moduli space of $\mathbb{C}^{*}$-Higgs bundles over $A$, is
a hyperKähler manifold. In particular, it has a holomorphic symplectic
form. Hence, both $\mathcal{M}^{H}$ and $X\,:=\,\mathcal{M}^{H}\otimes_{\mathbb{Z}}\Lambda_{T}$
are smooth symplectic varieties. Now, the Weyl group $W$ acts by
automorphisms of $X$ preserving the symplectic form, so that this
form descends to a symplectic form $\omega$ on the smooth part 
\[
\mathcal{M}_{A}^{H}(G)^{sm}\cong(X/W)^{sm}=(X\setminus\bigcup_{g\in W\setminus\{1\}}X^{g})/W,
\]
where $X^{g}$ is the fixed point locus for the action of the non-trivial
element $g\,\in\,W\setminus\{1\}$. Since $W$ acts by reflections
on $\Lambda_{T}$, the codimension of each $X^{g}$ coincides with
the dimension of $\mathcal{M}^{H}$, which is $2d$. Hence, the complement
of $(X/W)^{sm}$ in $X/W\,=\,\mathcal{M}_{A}^{H}(G)$ has also codimension
$2d\,\geq\,2$. Consequently, by \cite[Prop 2.4]{Be2}, $\omega$
extends to a symplectic form on $X/W$. The case of $\mathcal{M}_{A}^{R}(G)$
is analogous. The rational Gorenstein property is immediate from \cite[Prop. 1.3]{Be2}. 
\end{proof}
\begin{rem}
In the case $G\,=\,\text{SL}_{n}(\mathbb{C})$, the character variety
$\mathcal{M}_{A}^{R}(G)$ has been previously shown to be symplectic,
when $A$ is an elliptic curve (see \cite{BS}). 
\end{rem}

\subsection{Hilbert-Chow morphism and symplectic resolutions}

Given an algebraic variety $X$, the Hilbert-Chow morphism is: 
\[
\hilb^{n}(X)\,\longrightarrow\,\sym{}^{n}(X),
\]
where $\hilb^{n}(X)$ denotes the Hilbert scheme of zero dimensional
subschemes of length $n$ in the variety $X$, and the morphism sends
a scheme to its support. When $X$ is a smooth surface, the Hilbert-Chow
morphism is a desingularization of the symmetric product. From now
on, we only deal with the case $\GL\,=\,\GL(\mathbb{C})$.
\begin{cor}
\label{cor:resolution-E}When $C$ is an elliptic curve, 
\[
\hilb^{n}(T^{*}C)\,\longrightarrow\,\sym^{n}(T^{*}C)
\]
is a symplectic (crepant) resolution of singularities of the moduli
space $\mathcal{M}_{C}^{H}(\GL)$. In fact, this is the unique such
resolution. 
\end{cor}

\begin{proof}
This result appears in \cite[Thm. 4.9]{Ti}; it follows from Beauville
(\cite[Prop. 5]{Be1}), since $T^{*}C$ is an algebraic surface with
the natural Liouville symplectic form, and has trivial canonical bundle
(being a local statement, \cite[Prop. 5]{Be1} applies to non-projective
surfaces). The uniqueness follows from Fu-Namikawa's theorem \cite[Cor. 2.3]{FN}. 
\end{proof}
\begin{thm}
\label{thm:resolutions} Let $A$ be an abelian variety of dimension
$d\geq1$, and $n=2$ or $3$. Then 
\[
\hilb^{n}(T^{*}A^{\vee})\quad\quad(\text{respectively,}\quad\hilb^{n}(\mathbb{G}_{m}^{2d}))
\]
are resolutions of singularities of $\mathcal{M}_{A}^{H}(\GL)$ (respectively,
$\mathcal{M}_{A}^{R}(\GL)$). These are symplectic resolutions if
and only if $d\,=\,1$. 
\end{thm}

\begin{proof}
In \cite[Thm. 3.2.2]{C2}, Cheah showed that $\hilb^{2}(X)$ and $\hilb^{3}(X)$
are smooth varieties, when $X$ is smooth. Since $\mathcal{M}^{H}\,\cong\,T^{*}A^{\vee}$
and $\mathbb{G}_{m}^{2d}$ are smooth, the result is clear. See also
\cite[Rem. 5.11]{FT}. When $d\,=\,1$ the resolutions are symplectic
by Corollary \ref{cor:resolution-E}.

When $d\,>\,1$, we start by considering a general smooth symplectic
variety $X$ with an action by a finite group $F$ that preserves
the symplectic form. Then $X/F$ is $\mathbb{Q}$-\emph{factorial}
and normal. Hence, any component $E$ of the exceptional locus of
a proper resolution $\pi\,:\,Z\,\longrightarrow\,X/F$ is of codimension
1. On the other hand, if $\pi$ is a symplectic resolution, then by
the semi-smallness property we have $2\,=\,2\codim(E)\,\geq\,\codim(\pi(E))$
(see \cite[Lemma 2. 7]{Ka}). As $X/F$ is a normal variety, its singular
locus is of codimension $\,\geq\,2$, and therefore $\codim(\pi(E))\,=\,2$.
However, the singular locus of $X/F$ is contained in $p(Y)$, where
$p\,:\,X\,\longrightarrow\,X/F$ is the natural projection, 
\[
Y\,:=\,\bigcup_{g\in G\setminus\{1\}}X^{g},
\]
and $X^{g}\subset X$ is the subvariety fixed by the action of the
non-trivial element $g\in F\setminus\{1\}$. Hence $\codim(Y)\,\leq\,2$
and since all components of $Y$ are even dimensional (they are fixed
loci under a symplectic action), $Y$ has to be of pure codimension
2.

On the other hand, in view of the isomorphism 
\[
\mathcal{M}_{A}^{H}(\GL)\,\cong\,(\mathcal{M}^{H}\otimes_{\mathbb{Z}}\Lambda_{T})/W
\]
from Theorem \ref{thm:M^H2}, since $W$ acts by reflections on $\Lambda_{T}$,
the codimension of the fixed subvariety coincides with the dimension
of $\mathcal{M}^{H}$, which is $2d\,>\,2$. So, taking $X\,=\,\mathcal{M}^{H}\otimes_{\mathbb{Z}}\Lambda_{T}$
and $F\,=\,W$ in the previous paragraph we obtain a contradiction.
This proves that $\hilb^{n}(T^{*}A^{\vee})$ is not a symplectic resolution
of $\mathcal{M}_{A}^{H}(\GL)$ when $d\,>\,1$.

Analogously, $\hilb^{n}(\mathbb{G}_{m}^{2d})\,\longrightarrow\,\mathcal{M}_{A}^{R}(\GL)$
is not a symplectic resolution, for $d\,>\,1$. 
\end{proof}

\section{Topology of Resolutions via Hilbert Schemes}

In this section, we consider $G\,=\,{\rm GL}_{n}(\mathbb{C})$, and
compute the Poincaré polynomials of the natural desingularizations
of $\mathcal{M}_{A}^{H}(G)$ and of $\mathcal{M}_{A}^{R}(G)$ in all
cases when these are given by the Hilbert-Chow morphism. %

Our strategy is based on Cheah's identity for the generating series
of Serre polynomials of Hilbert schemes. This is enough to get the
Poincaré of the resolutions $\hilb^{n}(T^{*}A^{\vee})$ and $\hilb^{n}(\mathbb{G}_{m}^{2d})$,
given that the Hodge numbers on these moduli spaces have special symmetries.

\subsection{The Serre polynomial and Hilbert $E$-series}

The Serre polynomial (also known as $E$-polynomial) of a quasi-projective
complex variety $X$ is defined by: 
\[
E_{X}(u,v):=\mu_{X}^{c}(-1,u,v)=\sum_{k,p,q\geq0}(-1)^{k}\,h_{c}^{k,p,q}(X)\,u^{p}\,v^{q},
\]
where $h_{c}^{k,p,q}(X)\,=\,\dim_{\mathbb{C}}H_{c}^{k,p,q}(X,\,\mathbb{Q})$
and we use the super/subscript $c$ to refer to \emph{compactly supported
cohomology}.

Assuming that $X$ is smooth (but not necessarily projective), the
Poincaré polynomial $P_{X}(t)\,:=\,\sum_{k=0}^{2d}\dim_{\mathbb{C}}H^{k}(X,\,\mathbb{Q})\,t^{k}$
(for usual cohomology) can be deduced from $E_{X}$ in the pure cohomology
case. Recall that $X$ is said to have a \emph{pure Hodge structure}
if $h^{k,p,q}(X)=0$ unless $k=p+q$ (which, for $X$ smooth, is equivalent
to requiring $h_{c}^{k,p,q}(X)\,=\,0$ unless $k\,=\,p+q$). 
\begin{lem}
\label{lem:E-P}Let $E_{X}(u,\,v)\,\,\in\,\,\mathbb{Z}[u,\,v]$ be
the $E$-polynomial of a smooth quasi-projective variety $X$, of
complex dimension $d$. Then its Poincaré polynomial is given by 
\[
P_{X}(t)\,=\,t^{2d}E_{X}(-t^{-1},\,-t^{-1})
\]
whenever the cohomology of $X$ is pure. Moreover, $P_{X}(t)\,=\,E_{X}(-t,\,-t)$
if $X$ is projective. 
\end{lem}

\begin{proof}
This is a simple calculation. From the definition: 
\[
E_{X}(-t,\,-t)\,=\,\sum_{p,q=0}^{d}(-1)^{p+q}h_{c}^{p+q,p,q}(X)\,(-t)^{p}(-t)^{q}\,=\,\sum_{p,q=0}^{d}h_{c}^{p+q,p,q}(X)\,t^{p+q},
\]
which is exactly the compactly supported Poincaré polynomial $P_{X}^{c}(t)$
since $\sum_{p+q=k}h_{c}^{p+q,p,q}(X)=\dim H_{c}^{k}(X;\,\mathbb{Q})$.
Hence, the Lemma follows from Poincaré duality which implies $P_{X}^{c}(t)=t^{2d}P_{X}(t^{-1})$,
in the quasi-projective case, and $P_{X}=P_{X}^{c}$ when $X$ is
projective. 
\end{proof}
Let us now define the \emph{Hilbert $E$-series of $X$, }which collects
the Serre $E$-polynomials of $\hilb^{n}X$, for all $n\in\mathbb{N}$. 
\begin{defn}
Let $X$ be a complex quasi-projective variety. The \emph{Hilbert
$E$-series of $X$} is the formal power series (in the formal variable
$y$) defined by: 
\begin{equation}
\HE_{X}(u,v,y):=\sum_{n\geq0}E_{\hilb^{n}X}(u,v)\,y^{n}\in\mathbb{Z}[u,v]\llbracket y\rrbracket.\label{eq:HE-X}
\end{equation}
\end{defn}

In \cite{C1}, J. Cheah introduced a method for determining $\HE_{X}$
from the knowledge of $E_{X}$ and of $\HE_{(\mathbb{C}^{d},0)}$,
where $\hilb^{n}(\mathbb{C}^{d},0)$ are the \emph{punctual Hilbert
schemes} of subschemes of $\mathbb{C}^{d}$ which have length $n$,
dimension 0, and are supported at the origin. To write this result
in a convenient form, we introduce the so called \emph{plethystic
exponential }and\emph{ logarithm} functions. Given a power series
of the form $f(u,v,y)=\sum_{n\geq1}\sum_{p,q\geq0}\,c_{p,q,n}u^{p}v^{q}\,y^{n}\in\mathbb{Z}[u,v]\llbracket y\rrbracket$,
we define: 
\begin{equation}
\Pexp[f(u,v,y)]:=\exp\left(\sum_{n\geq1}\frac{1}{n}f(u^{n},v^{n},y^{n})\right)=\prod_{n\geq1}\prod_{p,q\geq1}(1-u^{p}v^{q}\,y^{n})^{-c_{p,q,n}}.\label{eq:Pexp}
\end{equation}
The equivalence of the two above expressions of $\Pexp$ follows from
the series expansion of the usual logarithm. The plethystic logarithm
$\Plog$ is the inverse function of $\Pexp$ (see for example \cite{FNZ1}). 
\begin{thm}[{{{{\cite[Main Thm. §2]{C1}}}}}]
\label{thm:Cheah} Let $X$ be smooth and quasi-projective of dimension
$d$. Then, we have: 
\[
\HE_{X}(u,\,v,\,y)\,=\,\Pexp\left(E_{X}(u,\,v)\,\Plog(\HE_{(\mathbb{A}^{d},0)}(u,\,v,\,y))\right).
\]
\end{thm}

\begin{proof}
This formula is due to Cheah \cite[S~2]{C1}, written in the form:
\begin{equation}
\HE_{X}(u,\,v,\,y)\,=\,\exp\left(\sum_{n\geq1}\frac{1}{n}E_{X}(u^{n},v^{n})\,g_{d}(u^{n},v^{n},y^{n})\right)\label{eq:Cheah}
\end{equation}
where $g_{d}(u,\,v,\,y)\,=\,\sum_{n\geq1}\sum_{p,q\geq0}\ a_{d}(p,\,q,\,n)\,u^{p}v^{q}y^{n}$,
and the coefficients $a_{d}(p,\,q,\,n)\in\mathbb{Q}$ are defined
by the infinite product expansion: 
\begin{equation}
\HE_{(\mathbb{C}^{d},0)}(u,\,v,\,y)\,=\,\prod_{n\geq1}\ \prod_{p,q\geq0}(1-u^{p}v^{q}\,y^{n})^{-a_{d}(p,\,q,\,n)}.\label{eq:def-a_d}
\end{equation}
Using the product form of \eqref{eq:Pexp}, Equation \eqref{eq:def-a_d}
is seen to be equivalent to 
\[
g_{d}(u,\,v,\,n)\,=\,\Plog(\HE_{(\mathbb{C}^{d},0)}(u,\,v,\,y)),
\]
and using the other form of \eqref{eq:Pexp} we get \eqref{eq:Cheah}. 
\end{proof}
To apply Theorem \ref{thm:Cheah} to a surface $S$, since $\dim S=\,2$,
we now compute $\Plog(\HE_{(\mathbb{C}^{2},\,0)})$. 
\begin{cor}
\label{cor:HE-S} Let $S$ be a quasi-projective smooth surface. Then:
\[
\HE_{S}(u,\,v,\,y)\,=\,\Pexp\left(\frac{E_{S}(u,\,v)\,y}{1-uvy}\right).
\]
\end{cor}

\begin{proof}
In \cite[Cor. (4) of Thm. 3.3.3]{C2} the following expansion is obtained:
\[
\sum_{n\geq0}P_{\hilb^{n}(\mathbb{C}^{2},\,0)}(t)\,y^{n}\,=\,\prod_{n\geq1}(1-t^{2n-2}y^{n})^{-1},
\]
so that $P_{\hilb^{n}(\mathbb{C}^{2},\,0)}(t)$ are polynomials in
$t^{2}$. Since $\hilb^{n}(\mathbb{C}^{2},\,0)$ is smooth and projective,
and has a cellular decomposition (\cite{C2}), its Hodge structure
is pure and balanced (i.e., $h_{c}^{k,p,q}\,=\,0$ if $p\,\neq\,q$,
or $k\,\neq\,p+q$) and Lemma \ref{lem:E-P} gives: 
\[
\HE_{(\mathbb{C}^{2},0)}(u,v,y)=\sum_{n\geq0}E_{\hilb^{n}(\mathbb{C}^{2},0)}(u,v)\,y^{n}=\prod_{n\geq1}(1-(uv)^{n-1}y^{n})^{-1}.
\]
Hence, by definition of $\Plog$: 
\[
\Plog\left(\HE_{(\mathbb{C}^{2},0)}(u,v,y)\right)=\sum_{n\geq1}\,(uv)^{n-1}y^{n}=\frac{y}{1-uvy}
\]
and the Corollary follows from Theorem \ref{thm:Cheah}. 
\end{proof}

\subsection{The case of elliptic curves}

Let $C$ be an elliptic curve. We now use Corollary \ref{cor:HE-S}
to obtain the Poincaré polynomials of $\hilb^{n}(T^{*}C)$, which
are crepant resolutions of $\mathcal{M}_{C}^{H}(\GL)$, by Corollary
\ref{cor:resolution-E}, since $C\,=\,C^{\vee}$. For that, note that
Equation \eqref{eq:HE-X}, for a smooth surface $S$ with pure cohomology,
implies: 
\[
\HE_{S}(-t^{-1},-t^{-1},yt^{4}):=\sum_{n\geq0}E_{\hilb^{n}S}(-t^{-1},-t^{-1})\,t^{4n}y^{n}=\sum_{n\geq0}P_{\hilb^{n}S}(t)\,y^{n}.
\]
Now, let $\mathcal{P}_{n}$ be the set of partitions of $n\in\mathbb{N}$,
and denote a partition as: 
\[
\lambda=[1^{\lambda_{1}}2^{\lambda_{2}}\cdots n^{\lambda_{n}}]\in\mathcal{P}_{n},
\]
where $\lambda_{j}\,\geq\,0$ is the number of parts of size $j\,=\,1,\,\cdots,\,n$,
so that $n\,=\,\sum_{j=1}^{n}j\lambda_{j}$. 
\begin{thm}
\label{thm:Poincare}Let $C$ be an elliptic curve and $n\in\mathbb{N}$.
The Poincaré polynomials of $\hilb^{n}(T^{*}C)$, the unique crepant
resolution of $\mathcal{M}_{C}^{H}(GL_{n}\mathbb{C})$, are given
by: 
\[
P_{\hilb^{n}(T^{*}C)}(t)=\frac{(t+1)^{2n}}{(t^{2}-1)^{n}}\sum_{\lambda\in\mathcal{P}_{n}}(t^{2\lambda_{1}}-1)\cdots(t^{2\lambda_{n}}-1)\,t^{2n-2|\lambda|},
\]
where $|\lambda|=\sum_{j=1}^{n}\lambda_{j}$ is the number of parts
of $\lambda\in\mathcal{P}_{n}$. 
\end{thm}

\begin{proof}
The surface $S=T^{*}C$ is homotopic to $C$, which gives (note the
factor $uv$ from compactly supported cohomology) $E_{T^{*}C}(u,v)=uv(u-1)(v-1)$.
So, from Corollary \ref{cor:HE-S}: 
\begin{eqnarray*}
\HE_{T^{*}C}(u,v,y) & = & \Pexp\left(\frac{(u-1)(v-1)\,uvy}{1-uvy}\right)\\
 & = & \frac{\Pexp(\frac{uvy}{1-uvy})\,\Pexp(\frac{(uv)^{2}y}{1-uvy})}{\Pexp(\frac{uv^{2}y}{1-uvy})\,\Pexp(\frac{u^{2}vy}{1-uvy})}\\
 & = & \prod_{n\geq1}\frac{(1-y^{n}u^{n}v^{n+1})(1-y^{n}u^{n+1}v^{n})}{(1-y^{n}u^{n+1}v^{n+1})(1-y^{n}u^{n}v^{n})}
\end{eqnarray*}
and 
\begin{equation}
\sum_{n\geq0}P_{\hilb^{n}T^{*}C}(t)\,y^{n}=\HE_{T^{*}C}(-t^{-1},-t^{-1},yt^{4})=\prod_{n\geq1}\frac{(1+y^{n}t^{2n-1})^{2}}{(1-y^{n}t^{2n-2})(1-y^{n}t^{2n})}.\label{eq:PHilb}
\end{equation}
To extract the coefficient of $y^{n}$ we use the (partial fraction)
expansion: 
\[
\frac{(1+\frac{z}{t})^{2}}{(1-\frac{z}{t^{2}})(1-z)}=1+\sum_{k\geq1}a_{k}(t)\,z^{k},
\]
where $a_{k}(t)=\frac{t+1}{t-1}(1-\frac{1}{t^{2k}})$, so that, replacing
$z$ successively by $(yt^{2})^{n}$, \eqref{eq:PHilb} becomes: 
\[
\frac{(1+yt)^{2}}{(1-y)(1-yt^{2})}\cdot\frac{(1+y^{2}t^{3})^{2}}{(1-y^{2}t^{2})(1-y^{2}t^{4})}\cdots=\sum_{n\geq0}(\frac{t+1}{t-1})^{n}\sum_{\lambda\in\mathcal{P}_{n}}(t^{2\lambda_{1}}-1)\cdots(t^{2\lambda_{n}}-1)\,t^{2n-2|\lambda|}y^{n},
\]
as wanted. 
\end{proof}
\begin{cor}
The Poincaré polynomial of $\hilb^{n}(\mathbb{G}_{m}^{2})$ coincides
with the one of $\hilb^{n}(T^{*}C)$, for all $n\ge1$. 
\end{cor}

\begin{proof}
The computations for $\hilb^{n}(\mathbb{G}_{m}^{2})$ follow the same
procedure as for $\hilb^{n}(T^{*}C)$, using now $E_{\mathbb{G}_{m}^{2}}(u,v)=(uv-1)^{2}$
and the fact that $\mathbb{G}_{m}^{2}=(\mathbb{C}^{*})^{2}$ has a
\emph{round} Hodge structure (its Hodge numbers verify $h^{k,p,q}\neq0$
only when $k=p=q$, see \cite{FS}). 
\end{proof}
\begin{rem}
Our computations were done in the quasi-projective setting; in the
projective context this agrees with results of Göttsche on Poincaré
polynomials of Hilbert schemes of smooth projective surfaces, and
also with the virtual motive computations in \cite{G}. Note also
that similar formulas for the case $G=\text{SL}_{n}(\mathbb{C})$
can be obtained by the same method. 
\end{rem}

\subsection{Higher dimensional abelian varieties}

We now turn to the computation of Poincaré polynomials of $\hilb^{n}(\mathcal{M}^{H})$
(respectively, $\hilb^{n}(\mathbb{G}_{m}^{2d})$) for abelian varieties
$A$ of higher dimension $d$, with $n\,=\,2$ and $3$ which, by
Theorem \ref{thm:resolutions} are natural resolutions of $\mathcal{M}_{A}^{H}(\GL)$
(respectively, $\mathcal{M}_{A}^{R}(\GL)$).

We use the following notation: $f_{1}(u,v)=E_{\mathcal{M}^{H}}(u,v)$,
$f_{1}^{*}(u,v):=f_{1}(u^{2},v^{2})$ and $f_{1}^{**}(u,v):=f_{1}(u^{3},v^{3})$.
Also, define, with $x=uv$, 
\begin{eqnarray}
f_{2} & := & E_{X}\,(E_{\hilb^{2}(\mathbb{A}^{2},0)}-1)=E_{X}\,(x+\cdots+x^{d-1})\nonumber \\
f_{3} & := & E_{X}\,(E_{\hilb^{3}(\mathbb{A}^{2},0)}-E_{\hilb^{2}(\mathbb{A}^{2},0)})=E_{X}\,\left(\sum_{1\leq i\leq j\leq d-1}x^{i+j}\right)\label{eq:f_i}
\end{eqnarray}
(see \cite{C2}). 
\begin{prop}
\label{prop:E-Hilb}The $E$-polynomials of $\hilb^{2}(\mathcal{M}^{H})$
and $\hilb^{3}(\mathcal{M}^{H})$ are given by: 
\begin{eqnarray*}
E_{\hilb^{2}(\mathcal{M}^{H})} & = & f_{2}+\frac{1}{2}(f_{1}^{2}+f_{1}^{*})\\
E_{\hilb^{3}(\mathcal{M}^{H})} & = & f_{3}+f_{1}f_{2}+\frac{f_{1}^{3}}{6}+\frac{f_{1}f_{1}^{*}}{2}+\frac{f_{1}^{**}}{3}.
\end{eqnarray*}
\end{prop}

\begin{proof}
We will use Theorem \ref{thm:Cheah}. It is not difficult to show
that 
\[
\Plog(\HE_{(\mathbb{C}^{d},0)}(u,v,y))=y+y^{2}(E_{\hilb^{2}(\mathbb{C}^{2},0)}-1)+y^{3}(E_{\hilb^{3}(\mathbb{C}^{2},0)}-E_{\hilb^{2}(\mathbb{C}^{2},0)})+\cdots
\]
so that multiplying by $E_{X}(u,v)$, we can write: 
\[
E_{X}(u,v)\Plog(\HE_{(\mathbb{C}^{d},0)}(u,v,y))=yf_{1}(u,v)+y^{2}f_{2}(u,v)+y^{3}f_{3}(u,v)+\cdots,
\]
where $f_{1}(u,v)=E_{X}(u,v)$. Hence we have 
\begin{eqnarray*}
\HE_{X}(u,v,y) & = & \Pexp(E_{X}(u,v)\,\Plog(\HE_{(\mathbb{C}^{d},0)}(u,v,y)))\\
 & = & \Pexp(yf_{1}(u,v)+y^{2}f_{2}(u,v)+y^{3}f_{3}(u,v)+\cdots)\\
 & = & \exp(yf_{1}+\frac{y^{2}}{2}f_{1}^{*}+\frac{y^{3}}{3}f_{1}^{**}+y^{2}f_{2}+y^{3}f_{3}+...)\\
 & = & 1+f_{1}y+(\frac{f_{1}^{2}}{2}+\frac{f_{1}^{*}}{2}+f_{2})y^{2}+(\frac{f_{1}^{3}}{6}+\frac{f_{1}f_{1}^{*}}{2}+\frac{f_{1}^{**}}{3}+f_{1}f_{2}+f_{3})y^{3}+\cdots
\end{eqnarray*}
which proves the Proposition. 
\end{proof}
Finally, the Poincaré polynomials of $\mathcal{M}_{A}^{H}(\GL)$ (respectively,
$\mathcal{M}_{A^{\vee}}^{R}(\GL)$), are obtained in three steps:
firstly, use $E_{T^{*}A}(u,v)=x^{d}(u-1)^{d}(v-1)^{d}$ and $E_{\mathbb{G}_{m}^{2d}}(u,v)=(x-1)^{2d}$
with $x=uv$, in the Equations \eqref{eq:f_i}. Secondly, get $E_{\hilb^{n}(\mathcal{M}^{H})}$
and $E_{\hilb^{n}(\mathcal{M}^{R})}$ for $n=2,3$ from Proposition
\ref{prop:E-Hilb}. Finally, since these are smooth varieties, use
Lemma \ref{lem:E-P} (or the corresponding relation between $E$ and
$P$, in case of round cohomology) to get $P_{\hilb^{n}(\mathcal{M}^{H})}(t)$
and $P_{\hilb^{n}(\mathcal{M}^{R})}(t)$ for $n=2,3$. 
\begin{rem}
The same method, using $E_{\mathbb{G}_{m}^{r}}(u,v)=(uv-1)^{r}$,
provides also the Poincaré polynomials of the natural resolutions
of the character varieties $\hom(\mathbb{Z}^{r},\text{GL}_{n}(\mathbb{C}))\quot\text{GL}_{n}(\mathbb{C})$
for odd $r$, even though, in this case, $\mathbb{Z}^{r}$ is not
a Kähler group, so it does not appear in the context of the non-abelian
Hodge correspondence. The analogous results of Theorem \ref{thm:Poincare}
and Proposition \ref{prop:E-Hilb}, for the case $G=\text{SL}_{n}(\mathbb{C})$
are also easy to obtain with the same methods. 
\end{rem}

\end{document}